\documentclass{article}

\usepackage{amssymb, amsmath, amsthm, enumerate, tikz}

\theoremstyle{plain}
\newtheorem{theorem}{Theorem}[section]
\newtheorem{lemma}[theorem]{Lemma}
\newtheorem{proposition}[theorem]{Proposition}

\theoremstyle{definition}
\newtheorem{definition}[theorem]{Definition}

\newenvironment{problem}[3]{
\begin{quote}
#1 \nopagebreak \\
Instance: #2\\
Question: #3}
{\end{quote}}

\newcommand{\iso}{\cong}
\newcommand{\del}{\backslash}
\newcommand{\cl}{\operatorname{cl}}
\newcommand{\dash}{\nobreakdash-\hspace{0pt}}
\newcommand{\npc}{\textup{NP}\nobreakdash-complete}

\newcommand{\mcal}[1]{\ensuremath{\mathcal{#1}}}
\newcommand{\tsc}[1]{\textsc{#1}}
\newcommand{\tup}[1]{\textup{#1}}

\title{Matroid complexity and non-succinct descriptions}

\author{Dillon Mayhew\\
Mathematical Institute\\
University of Oxford\\
St Giles\\
OX1 3LB\\
Oxford, U.K.
\thanks{Current address: School of Mathematics,
Victoria University of Wellington,
P.O. BOX 600,
Wellington 6140,
New Zealand.
(email {\tt dillon.mayhew@mcs.vuw.ac.nz})}}

\begin{document}

\maketitle

\begin{abstract}
We investigate an approach to matroid complexity that
involves describing a matroid via a list of independent
sets, bases, circuits, or some other family of subsets
of the ground set.
The computational complexity of algorithmic problems under this
scheme appears to be highly dependent on the choice of input-type.
We define an order on the various methods of description, and we show how
this order acts upon ten types of input.
We also show that under this approach several natural algorithmic problems
are complete in classes thought not to be equal to \tup{P}.
\end{abstract}

\section{Introduction}
\label{sct1}

The study of matroid-theoretical algorithmic problems and their
complexity has been dominated by two approaches.
The first approach was implicitly used by Edmonds~\cite{edmonds3} and
later developed by Hausmann and Korte~\cite{hausmann2,hausmann3}.
It uses the idea of an ordinary Turing machine augmented with an oracle.
Suppose that the subject of the computation is a matroid on the
ground set $E$.
When queried about a subset, $X$, of $E$, the oracle
returns in unit time some information about $X$.
That information might be the rank of $X$, or an answer to
the question ``Is $X$  independent?'', to mention the two
most widely used oracles.

The other approach to matroid complexity uses the standard model
of a Turing machine, but considers as its input only a restricted
class of matroids that can be represented by some `succinct' structure,
for instance, a graph, or a matrix over a field.

A third approach to the study of matroid complexity would, in some
ways, be more natural.
A matroid is essentially a finite set with a structured family of
subsets (we shall not consider infinite matroids).
An obvious way to describe a matroid to a Turing machine is to
simply list the subsets that belong to this family.
The advantage of this approach is that any matroid could be
received as input, and not merely those matroids which belong
to some restricted class.

The reason this approach has not received as much attention is
the concern that the input will be too large.
A result of Knuth's~\cite{knuth} says that if $f(n)$ is the number of
non-isomorphic matroids on a set of size $n$, then there is a
constant $c$ such that
\begin{displaymath}
\log_{2}\log_{2} f(n) \geq n - \frac{3}{2}\log_{2} n + c \log_{2}\log_{2} n
\end{displaymath}
for sufficiently large values of $n$.
It follows that if $\Sigma$ is a finite alphabet and
$\sigma: \mcal{M} \rightarrow \Sigma^{*}$ is an injective encoding
scheme that takes the set of all matroids to words in $\Sigma$,
then there can be no polynomial function $p$ such that
$|\sigma(M)|$, the length of the word $\sigma(M)$, is bounded above by
$p(|E(M)|)$ for every matroid $M$.
Thus if we wish to use this type of encoding function we must in some
sense abandon the cardinality of the ground set of a matroid as a
measure of its `size'.
Historically the concern with this type of scheme has been that,
because of the large size of the input, the class of problems that
can be solved in polynomial time will be artificially inflated, and that
therefore, in a trivial way, all algorithmic problems for matroids will
be tractable.

In this paper we show that the situation is apparently more subtle
than this.
We show that several natural matroid problems are complete in classes
thought not to be equal to \tup{P}, even using an encoding scheme that
works for all matroids: for instance, the scheme that describes a
matroid by listing its independent sets.

A quirk of this approach to matroid complexity is that the
difficulty of a computational problem appears to vary widely
according to the type of input. A problem may be solvable in
polynomial time if the input takes the form of a list of bases, but
if the input is a list of circuits the same problem may be \npc.

Before we examine the complexity of matroid-theoretical problems,
we define an order on types of input, and we show how this order
acts upon a set of ten natural methods of description. This work
is an analogue of that done by Hausmann and Korte~\cite{hausmann},
and Robinson and Welsh~\cite{robinson}, comparing different types
of oracles.

Our references for basic concepts, notation, and terminology will
be Oxley~\cite{oxley} with regards to matroids, and Garey
and Johnson~\cite{garey} with regards to complexity theory.

\section{Various types of inputs}
\label{sct2}

The ten types of input that we consider are as follows:
\tsc{Rank}, \tsc{Independent Sets}, \tsc{Spanning Sets},
\tsc{Bases}, \tsc{Flats}, \tsc{Circuits}, \tsc{Hyperplanes},
\tsc{Non-Spanning Circuits}, \tsc{Dependent Hyperplanes},
and \tsc{Cyclic Flats}.

Of these, some need little explanation. A list of the
independent sets, spanning sets, bases, flats, circuits, or
hyperplanes of a matroid uniquely specifies that matroid.
Thus the corresponding forms of input will consist simply
of lists of the appropriate subsets.

The \tsc{Rank} input will list each subset of the ground set,
along with its rank.
The \tsc{Non-Spanning Circuits} input for a
matroid $M$ will specify the rank of $M$ and list all its
non-spanning circuits.
Dually, the \tsc{Dependent Hyperplanes} input 
will specify the rank of $M$ and list its dependent
hyperplanes.

A \emph{cyclic flat} is a flat that is also a
(possibly empty) union of circuits.
It is known that listing the cyclic flats and specifying their
ranks completely determines a matroid.
Therefore the \tsc{Cyclic Flats} input will list each cyclic flat,
along with its rank.

(Note that we have not specified how to describe some exceptional
cases, such as a matroid $M$ with no non-spanning circuits.
In this particular case we will assume that the
\tsc{Non-Spanning Circuits} description lists only the rank of $M$.
Other exceptional cases are easily dealt with in a similar way.)

Suppose that $f$ and $g$ are two functions on the positive
integers.
If there exist constants, $c_{1}$ and $c_{2}$, and an
integer, $N$, such that $c_{1}g(n) \leq f(n) \leq c_{2}g(n)$ for
every positive integer $n \geq N$, then we shall write
$f = \Theta(g)$.
Equivalently, $f = \Theta(g)$ if and only if $f = O(g)$ and $g = O(f)$.

Suppose that $M$ is a matroid on a ground set of size $n$,
and that \tsc{Input} is one of the types of input discussed above.
Let $(M,\, \tsc{Input})$ be a word that describes $M$ via \tsc{Input}.
We shall assume that a subset of the ground set is specified by
its characteristic vector.
Thus if the \tsc{Input} description involves listing $i$ subsets of
$E(M)$ then $|(M,\, \tsc{Input})| \geq ni$.
We shall consider only `reasonable' encoding schemes
(i.e. those that do not allow, for instance, padding of words).
It follows easily that $|(M,\, \tsc{Input})| = \Theta(ni)$.

Obviously there are many other natural ways of describing a
matroid, but many are related in a fairly trivial way to one
of the methods we have already discussed.
For example, a matroid can be described by listing its cocircuits,
but the cocircuits are exactly the complements of the hyperplanes.

\section{A comparison of inputs}
\label{sct3}

It is natural to ask whether one form of description is
intrinsically more compact than another. In this section we
attempt to answer that question.

\begin{definition}
\label{COMdef1}
Suppose that \tsc{Input1} and \tsc{Input2} are two
methods for describing a matroid.
Then $\tsc{Input1} \leq \tsc{Input2}$ if there exists a
polynomial-time Turing machine which will produce
$(M,\, \tsc{Input2})$ given $(M,\, \tsc{Input1})$
for any matroid $M$.
\end{definition}

Suppose that $\tsc{Input1} \leq \tsc{Input2}$.
If a problem is in \tup{P} for descriptions via \tsc{Input2} then
clearly it is in \tup{P} for \tsc{Input1}.
Similarly, if a problem is \npc\ for \tsc{Input1} then the same
problem is \tup{NP}\dash hard for \tsc{Input2}.

It is clear that $\leq$ is both reflexive and transitive.
The rest of this section will be devoted to proving
Theorem~\ref{thm1}.

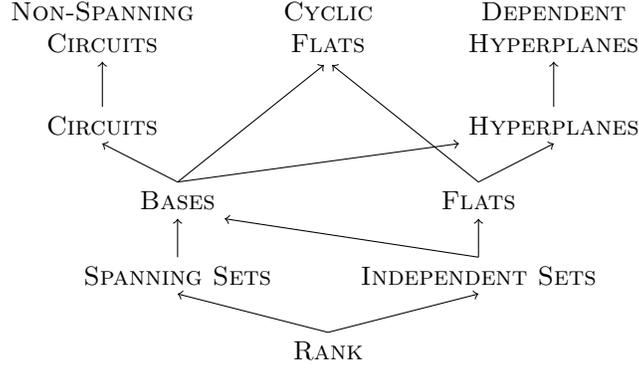
\begin{figure}[htb]

\centering

\begin{tikzpicture}
\draw (0,0) node(rank){\tsc{Rank}};
\draw (-2, 1) node(span){\tsc{Spanning Sets}};
\draw (2, 1) node(indp){\tsc{Independent Sets}};
\draw (-2, 2) node(base){\tsc{Bases}};
\draw (2, 2) node(flat){\tsc{Flats}};
\draw (-3, 3) node(circ){\tsc{Circuits}};
\draw (3, 3) node(hype){\tsc{Hyperplanes}};
\draw (-3, 4.3) node[text width=3cm, text centered](nspc)
{\tsc{Non-Spanning Circuits}};
\draw (3, 4.3) node[text width=3cm, text centered](deph)
{\tsc{Dependent Hyperplanes}};
\draw (0, 4.3) node[text width=1.2cm, text centered](cycf)
{\tsc{Cyclic Flats}};

\draw[->] (rank.north) -- (span.south);
\draw[->] (rank.north) -- (indp.south);
\draw[->] (span.north) -- (base.south);
\draw[->] (indp.north) -- (base.south east);
\draw[->] (indp.north) -- (flat.south);
\draw[->] (base.north) -- (circ.south);
\draw[->, shorten >=2pt] (base.north) -- (cycf.south);
\draw[->] (base.north) -- (hype.south west);
\draw[->, shorten >=2pt] (flat.north) -- (cycf.south);
\draw[->] (flat.north) -- (hype.south);
\draw[->] (circ.north) -- (nspc.south);
\draw[->] (hype.north) -- (deph.south);
\end{tikzpicture}

\caption{An ordering of input types.}

\label{COMfig2}

\end{figure}

\begin{theorem}
\label{thm1}
The ten types of input listed in Section~\ref{sct2} are
ordered by $\leq$ according to the Hasse diagram in
Figure~\ref{COMfig2}.
\end{theorem}

Theorem~\ref{thm1} will follow from the
subsequent lemmas.

\begin{lemma}
\label{lem1}
$\tsc{Rank} \leq \tsc{Spanning Sets}$ and
$\tsc{Rank} \leq \tsc{Independent Sets}$.
\end{lemma}

\begin{proof}
To find $(M,\, \tsc{Spanning Sets})$ given $(M,\, \tsc{Rank})$, a
Turing machine need only write onto its output those subsets whose
rank equals $r(M)$.
Similarly, to produce $(M,\, \tsc{Independent Sets})$ the machine
need only write those subsets of $E(M)$ whose rank is equal to their size.
\end{proof}

\begin{lemma}
\label{lem2}
$\tsc{Spanning Sets} \leq \tsc{Bases}$ and
$\tsc{Independent Sets} \leq \tsc{Bases}$.
\end{lemma}

\begin{proof}
The bases of a matroid are exactly the minimal spanning sets and 
the maximal independent sets, so the result follows easily.
\end{proof}

\begin{lemma}
\label{lem3}
$\tsc{Independent Sets} \leq \tsc{Flats}$.
\end{lemma}

\begin{proof}
Given the list of independent sets of $M$, we can create the
list of flats by finding, in turn, the closure of each independent
set, and then eliminating duplications.
To find the closure of the independent set $I$, we check each
set of the form $I \cup e$, where $e \notin I$, to see if it
is independent.
The element $e$ belongs to $\cl(I)$ if and only if $I \cup e$
is dependent.
Since $|(M,\, \tsc{Independent Sets})| = \Theta(ni)$, where
$|E(M)| = n$ and the number of independent sets is $i$,
it is easy to see that this entire procedure can be accomplished
in polynomial time.
\end{proof}

\begin{lemma}
\label{lem4}
$\tsc{Bases} \leq \tsc{Circuits}$.
\end{lemma}

\begin{proof}
If $B$ is a basis of a matroid $M$ and
$e \in E(M) - B$, then $B \cup e$ contains a unique
circuit $C(e,\, B)$, known as the
\emph{fundamental circuit of $e$ with respect to $B$}.
Since every circuit is a fundamental circuit with respect to
some basis, we can construct the list of circuits of $M$, given
the list of bases, by creating the list of fundamental circuits
with respect to each of the bases, and then eliminating duplications.

If $B$ is a basis of $M$, and $e \in E(M) - B$, we find $C(e,\, B)$ by
comparing $(B \cup e) - f$ against the list of bases for each $f \in B$.
The element $f$ is in $C(e,\, B)$ if and only if $(B \cup e) - f$ is a
basis.
\end{proof}

The next result follows easily using duality and
Lemma~\ref{lem4}.

\begin{lemma}
\label{lem5}
$\tsc{Bases} \leq \tsc{Hyperplanes}$.
\end{lemma}

Let the number of bases and cyclic flats in the matroid $M$
be denoted by $b(M)$ and $z(M)$, respectively.

\begin{proposition}
\label{prop4}
Let $M$ be a matroid. Then $z(M) \leq b(M)$.
\end{proposition}

\begin{proof}
The proof will be by induction on $|E(M)|$.
If $|E(M)| = 0$, then $b(M) = z(M)= 1$, so the proposition holds.

Let $M$ be a matroid such that $|E(M)| = n > 0$, and assume that
the proposition holds for all matroids on ground sets of $n-1$ elements.
We may assume that there exists an element $e \in E(M)$ such that
$e$ is neither a loop nor a coloop, for otherwise $b(M) = z(M) = 1$.

Let $b_{e}(M)$ be the number of bases of $M$ that contain $e$, and
let $b_{\bar{e}}(M)$ be the number of bases of $M$ that avoid $e$.
Then $b_{e}(M) = b(M/e)$ and $b_{\bar{e}}(M) = b(M \del e)$.
Similarly, let $z_{e}(M)$ be the number of cyclic flats of $M$ that
contain $e$, and let $z_{\bar{e}}(M)$ be the number of cyclic flats
that do not contain $e$.
It is easy to see that any cyclic flat that does not contain $e$
is also a cyclic flat of $M \del e$, so $z_{\bar{e}}(M) \leq z(M \del e)$.
Moreover, if $Z$ is a cyclic flat of $M$, and $e \in Z$, then $Z - e$
is a cyclic flat of $M / e$.
Thus $z_{e}(M) \leq z(M / e)$, and hence
\begin{displaymath}
z(M) = z_{e}(M) + z_{\bar{e}}(M) \leq z(M/e) + z(M \del e).
\end{displaymath}
By the inductive hypothesis, $z(M/e) \leq b(M/e)$ and
$z(M \del e) \leq b(M \del e)$, so
\begin{displaymath}
z(M) \leq b(M / e) + b(M \del e) = b_{e}(M) + b_{\bar{e}}(M) = b(M). \qedhere
\end{displaymath}
\end{proof}

\begin{lemma}
\label{lem14}
$\tsc{Bases} \leq \tsc{Cyclic Flats}$.
\end{lemma}

\begin{proof}
Our algorithm for generating the list of cyclic flats, given
the list of bases, will start by constructing the closures of
all circuits.
At each repetition of the loop the algorithm will
find the closure of the union of every pair of cyclic flats
already on the list.

Suppose that $M$ is a matroid on the ground set $E$.
Note that we can check in polynomial time whether a
set is independent by comparing it against the family of bases.
Hence if $A \subseteq E$ we can use the greedy algorithm to find a
basis $I$ of $A$.
The element $e \notin A$ is in $\cl(A)$ if and only if
$I \cup e$ is dependent.
It follows that we can find $\cl(A)$ in polynomial time.

From Lemma~\ref{lem4} we see that we can find
the list of circuits in polynomial time, given the list of bases.
Hence we can also construct the list of closures of circuits, add
the closure of the empty set, and then eliminate duplications from
the list in polynomial time.
This completes the preprocessing the algorithm will do before
entering the loop.

Suppose that $Z_{1},\ldots, Z_{t}$ is the
list of cyclic flats that has been constructed after
the loop has been repeated $i$ times.
In the next repetition of the loop the algorithm will take
each pair $\{Z_{j},\, Z_{k}\}$ of cyclic flats
and find $\cl(Z_{j} \cup Z_{k})$.
At the completion of the loop the algorithm will add these
new cyclic flats to the list and then eliminate duplications.
By Proposition~\ref{prop4} the length of the list
at the beginning of the loop will never exceed $b$, the number of bases.
Thus in each repetition of the loop the algorithm must find
the closure of at most $b^{2}$ unions of flats.
It is easy to see that after $r(M)$ repetitions of the
loop the algorithm will have found every cyclic flat.
Finding the rank of these flats can clearly be accomplished in
polynomial time, using the greedy algorithm.
Therefore the algorithm can construct $(M,\, \tsc{Cyclic Flats})$
in polynomial time.
\end{proof}

\begin{lemma}
\label{lem6}
$\tsc{Flats} \leq \tsc{Hyperplanes}$.
\end{lemma}

\begin{proof}
The flat $F$ is a hyperplane of $M$ if and only if there is no
flat, $F'$, such that $F \subset F' \subset E(M)$.
This clearly leads to a polynomial-time algorithm.
\end{proof}

\begin{lemma}
\label{lem7}
$\tsc{Flats} \leq \tsc{Cyclic Flats}$.
\end{lemma}

\begin{proof}
A flat $F$ fails to be a cylic flat if and only if there is
an element $e \in F$ such that $F - e$ is also a flat.
This clearly leads to a polynomial time algorithm.
\end{proof}

\begin{lemma}
\label{lem8}
$\tsc{Circuits} \leq \tsc{Non-Spanning Circuits}$.
\end{lemma}

\begin{proof}
Given $(M,\, \tsc{Circuits})$ a Turing machine can check in
polynomial time whether a subset $A \subseteq E(M)$ is independent.
Thus such a machine could find a basis of $M$ by using the
greedy algorithm.
Once the rank of $M$ is known the rest follows easily.
\end{proof}

The next lemma follows easily using duality and Lemma~\ref{lem8}.

\begin{lemma}
\label{lem9}
$\tsc{Hyperplanes} \leq \tsc{Dependent Hyperplanes}$.
\end{lemma}

To complete the proof of Theorem~\ref{thm1} we must
show that if \tsc{Input1} and \tsc{Input2} are two types
of input and the preceding results do not imply that
$\tsc{Input1} \leq \tsc{Input2}$, then
$\tsc{Input1} \nleq \tsc{Input2}$.
It is clear that if $\tsc{Input1} \leq \tsc{Input2}$ then there
must exist a polynomial $p$ such that
\begin{displaymath}
|(M,\, \tsc{Input2})| \leq p(|(M,\, \tsc{Input1})|)
\end{displaymath}
for any matroid $M$.
Thus if there exists a family of matroids $M_{1},\, M_{2},\, M_{3},\ldots$
such that $\{|(M_{i},\, \tsc{Input2})|\}_{i \geq 1}$ is not
polynomially bounded by the sequence
$\{|(M_{i},\, \tsc{Input1})|\}_{i \geq 1}$ then
$\tsc{Input1} \nleq \tsc{Input2}$.

The following proposition follows immediately from
transitivity.

\begin{proposition}
\label{prop2}
Suppose that $\tsc{Input1} \nleq \tsc{Input2}$.
If $\tsc{Input1} \leq \tsc{Input3}$, and
$\tsc{Input4} \leq \tsc{Input2}$, then
$\tsc{Input3} \nleq \tsc{Input4}$.
\end{proposition}

Using Proposition~\ref{prop2}, it is a relatively
simple matter to check that the proof of Theorem~\ref{thm1}.
will be completed by verifying the following cases.

\begin{enumerate}[1.]
\item\label{1} $\tsc{Spanning Sets} \nleq \tsc{Flats}$
\item\label{2} $\tsc{Independent Sets} \nleq \tsc{Spanning Sets}$
\item\label{3} $\tsc{Flats} \nleq \tsc{Non-Spanning Circuits}$
\item\label{4} $\tsc{Circuits} \nleq \tsc{Dependent Hyperplanes}$
\item\label{5} $\tsc{Hyperplanes} \nleq \tsc{Cyclic Flats}$
\item\label{6} $\tsc{Circuits} \nleq \tsc{Cyclic Flats}$
\item\label{7} $\tsc{Non-Spanning Circuits} \nleq \tsc{Circuits}$
\item\label{8} $\tsc{Dependent Hyperplanes} \nleq \tsc{Hyperplanes}$
\item\label{9} $\tsc{Cyclic Flats} \nleq \tsc{Dependent Hyperplanes}$
\end{enumerate}

\begin{lemma}[Case~\ref{1}]
\label{lem10}
$\tsc{Spanning Sets} \nleq \tsc{Flats}$.
\end{lemma}

\begin{proof}
For $n \geq 1$ let $M_{n}$ be isomorphic to $U_{n-1,n}$, the
$n$\dash element uniform matroid of rank $n - 1$.
The number of spanning sets of $M_{n}$ is $n + 1$, whereas the number of flats
is $2^{n} - n$.
Thus $|(M_{n},\, \tsc{Spanning Sets})| = \Theta(n^{2})$, while
$|(M_{n},\, \tsc{Flats})| = \Theta(n2^{n})$.
\end{proof}

\begin{lemma}[Case~\ref{2}]
\label{lem11}
$\tsc{Independent Sets} \nleq \tsc{Spanning Sets}$.
\end{lemma}

\begin{proof}
For $n \geq 1$ let $M_{n}$ be isomorphic to $U_{1,n}$.
The number of independent sets in $M_{n}$ is $n + 1$,
while the number of spanning sets is $2^{n} - 1$.
\end{proof}

We denote the truncation of the matroid $M$ by $T(M)$.
If $m$ is a positive integer then we define $mU_{r,n}$ to
be the matroid obtained by replacing each element in $U_{r,n}$
with a parallel class of size $m$.

\begin{lemma}[Case~\ref{3}]
\label{lem15}
$\tsc{Flats} \nleq \tsc{Non-Spanning Circuits}$.
\end{lemma}

\begin{proof}
For $n \geq 3$ define $M_{n}$ to be $T(nU_{n-1,n} \oplus U_{2,2})$.
Note that $M_{n}$ contains $n^{2} + 2$ elements.
There are $n + 2$ parallel classes in $M_{n}$.
It follows that the number of flats of $M_{n}$ is at most $2^{n + 2}$.
However the number of non-spanning circuits of $M_{n}$ is exactly
$n^{n} + n\binom{n}{2}$.
Thus $|(M_{n},\, \tsc{Flats})| = O(n^{2}2^{n+2})$ while
$|(M_{n},\, \tsc{Non-Spanning Circuits})| = \Theta(n^{n+2})$.
\end{proof}

\begin{lemma}[Case~\ref{4}]
\label{lem16}
$\tsc{Circuits} \nleq \tsc{Dependent Hyperplanes}$.
\end{lemma}

\begin{proof}
It is not difficult to see that a polynomial-time
algorithm that constructs $(M,\, \tsc{Dependent Hyperplanes})$
from $(M,\, \tsc{Circuits})$ for any matroid $M$ can
be used to show that
$\tsc{Hyperplanes} \leq \tsc{Non-Spanning Circuits}$.
This contradicts Proposition~\ref{prop2}, as
$\tsc{Flats} \leq \tsc{Hyperplanes}$ by
Lemma~\ref{lem6}, and
$\tsc{Flats} \nleq \tsc{Non-Spanning Circuits}$
by Lemma~\ref{lem15}.
\end{proof}

\begin{lemma}[Case~\ref{5}]
\label{lem18}
$\tsc{Hyperplanes} \nleq \tsc{Cyclic Flats}$.
\end{lemma}

\begin{proof}
For $n \geq 3$, let $M_{n}$ be isomorphic to $2U_{n-1,n}$.
The hyperplanes of $M_{n}$ are exactly the sets of $n-2$
parallel classes, while any non-empty set of parallel classes
is a cyclic flat as long as it does not contain exactly
$n - 1$ such classes.
Thus the number of hyperplanes is $(n^{2}-n)/2$
while the number of cyclic flats is $2^{n} - n - 1$.
\end{proof}

\begin{lemma}[Case~\ref{6}]
\label{lem19}
$\tsc{Circuits} \nleq \tsc{Cyclic Flats}$.
\end{lemma}

\begin{proof}
Given $(M,\, \tsc{Hyperplanes})$ we can certainly find
$(M^{*},\, \tsc{Circuits})$ in polynomial time.
Also, given $(M^{*},\, \tsc{Cyclic Flats})$ we can
find the cyclic flats of $M$ in polynomial time,
since the cyclic flats of $M^{*}$ are the complements
of the cyclic flats of $M$.
Moreover it is easy to see that given
$(M^{*},\, \tsc{Circuits})$ we can find the rank of any
subset in $M$ in polynomial time.

Suppose that $\tsc{Circuits} \leq \tsc{Cyclic flats}$.
The above discussion implies that
$\tsc{Hyperplanes} \leq \tsc{Cyclic Flats}$, in
contradiction to Lemma~\ref{lem18}.
\end{proof}

\begin{lemma}[Case~\ref{7}]
\label{lem20}
$\tsc{Non-Spanning Circuits} \nleq \tsc{Circuits}$.
\end{lemma}

\begin{proof}
If $M_{n}$ is isomorphic to $U_{n,2n}$ then $M_{n}$ contains no
non-spanning circuits, so by definition
$|(M_{n},\, \tsc{Non-Spanning Circuits})| = O(n)$.
On the other hand, the number of circuits is
$\binom{2n}{n+1}$, which is exponential in $n$.
\end{proof}

The next lemma follows using duality and Lemma~\ref{lem20}.

\begin{lemma}[Case~\ref{8}]
\label{lem21}
$\tsc{Dependent Hyperplanes} \nleq \tsc{Hyperplanes}$.
\end{lemma}

\begin{lemma}[Case~\ref{9}]
\label{lem17}
$\tsc{Cyclic Flats} \nleq \tsc{Dependent Hyperplanes}$.
\end{lemma}

\begin{proof}
For $n \geq 2$ let $M_{n}$ be the matroid obtained by adding a
single parallel element to a member of the ground set of $U_{n,2n}$.
The only cyclic flats of $M_{n}$ are the empty set, the non-trivial
parallel class, and the entire ground set.
Thus $|(M_{n},\, \tsc{Cyclic Flats})| = \Theta(n)$.
However any hyperplane that contains the non-trivial parallel class is
dependent, so the number of such hyperplanes is $\binom{2n-1}{n-2}$.
\end{proof}

With this lemma we have completed the proof of
Theorem~\ref{thm1}.

\section{Matroid intersection}
\label{sct5}

It is easy to see that given $(M,\, \tsc{Input})$,
it is possible to determine in polynomial time whether or not
a subset of $E(M)$ is independent in $M$.
(Henceforth we assume \tsc{Input} to be one of the types of input
discussed in Section~\ref{sct2}.)
Hausmann and Korte~\cite{hausmann}, and Robinson and
Welsh~\cite{robinson} show that the standard matroid
oracles can be efficiently simulated by an independence
oracle.
It follows from these observations that if a computational
problem can be solved by an oracle Turing machine in time that
is bounded by $p(n)$ for any $n$\dash element matroid, where
$p$ is a fixed polynomial, then the same problem can be solved
in polynomial time by a Turing machine which receives
$(M,\, \tsc{Input})$ as its input.

The converse is not true.
Consider the problem of deciding whether a matroid is uniform.
Robinson and Welsh~\cite{robinson} note that a Turing machine
equipped with an oracle cannot solve this problem in time bounded
by a polynomial function of the size of the ground set.
In contrast, deciding whether $M$ is uniform given
$(M,\, \tsc{Input})$ is trivial.

However, there do exist matroid-theoretical problems which are
probably not solvable in polynomial time, even when the input consists
of a list of some family of subsets.
One of these is $3$\dash MATROID INTERSECTION.

\begin{problem}
{$3$\dash MATROID INTERSECTION}
{An integer $k$ and $(M_{i},\, \tsc{Input})$
for $1 \leq i \leq 3$, where $M_{1}$, $M_{2}$, and $M_{3}$
are matroids with a common ground set $E$.}
{Does there exist a set $A \subseteq E$ such
that $|A| = k$ and $A$ is independent in
$M_{1}$, $M_{2}$, and $M_{3}$?}
\end{problem}

The fact that this problem is \npc\ was first observed by
Lawler~\cite{lawler2}.
He does not specify a form of matroid input, but he remarks
that the problem is \npc\ for partition matroids
(direct sums of rank\dash one uniform matroids).
It is clear from his comments that a partition matroid is to be
described via the partition of its ground set into connected components.
We sketch a modified version of his proof here.

\begin{theorem}
\label{thm5}
If $\tsc{Circuits} \leq \tsc{Input}$ or if
$\tsc{Hyperplanes} \leq \tsc{Input}$, then
\tup{$3$\dash MATROID INTERSECTION} is \npc.
However, if $\tsc{Input} \leq \tsc{Bases}$ then
\tup{$3$\dash MATROID INTERSECTION} is in \tup{P}\@.
\end{theorem}

\begin{proof}
Obviously the problem is in \tup{NP}.
It suffices to prove \tup{NP}\dash completeness
only in the case that $\tsc{Input} = \tsc{Circuits}$ or
$\tsc{Input} = \tsc{Hyperplanes}$.
We provide a reduction from the following \npc\ problem.

\begin{problem}
{$3$\dash DIMENSIONAL MATCHING}
{A set of triples,
$M \subseteq X_{1} \times X_{2} \times X_{3}$, where
$X_{1}$, $X_{2}$, and $X_{3}$ are pairwise disjoint sets having
the same cardinality.}
{Does $M$ contain a \emph{matching}?
(A subset $M' \subseteq M$, such that every
element in $X_{1} \cup X_{2} \cup X_{3}$ is contained in exactly
one triple in $M'$.)}
\end{problem}

Let $M \subseteq X_{1} \times X_{2} \times X_{3}$ be an
instance of $3$\dash DIMENSIONAL MATCHING, and suppose
that $X_{i} = \{x_{1}^{i},\ldots, x_{s}^{i}\}$ for $1 \leq i \leq 3$.
Suppose that $M$ contains $t$ triples, $p_{1},\ldots, p_{t}$.
We construct three partition matroids, $M_{1}$, $M_{2}$, and $M_{3}$,
on the ground set $E = \{e_{1},\ldots, e_{t}\}$.
Each matroid $M_{i}$ contains $s$ connected components, corresponding to the
elements of $X_{i}$.
The connected component corresponding to $x_{j}^{i}$ is
equal to $\{e_{k} \mid p_{k}\ \mbox{contains}\ x_{j}^{i},\, 1 \leq k \leq t\}$.
It is clear that $M$ contains a matching if and only if
$M_{1}$, $M_{2}$, and $M_{3}$ contain a common independent set
of size $s$.

It remains to show that $(M_{i},\, \tsc{Circuits})$ and
$(M_{i},\, \tsc{Hyperplanes})$
can be constructed in polynomial time.
Since the number of circuits or hyperplanes in a partition matroid is
at most quadratic in the size of the ground set this is easily done.

If $\tsc{Input} = \tsc{Bases}$ then we can find a common
independent set of maximum size by considering the
intersection of every triple of bases from the three
matroids.
This can clearly be done in polynomial time, so
$3$\dash MATROID INTERSECTION is in \tup{P} if
$\tsc{Input} \leq \tsc{Bases}$.
\end{proof}

Theorem~\ref{thm5} shows that $3$\dash MATROID INTERSECTION
is either \npc\ or in \tup{P} for all but two of the methods of
input described in Section~\ref{sct2}:
The status of the problem is open for the case that
$\tsc{Input} = \tsc{Cyclic Flats}$ or
$\tsc{Input} = \tsc{Flats}$.

\section{The isomorphism problem}
\label{sct6}

The following computational problem has attracted much
attention.

\begin{problem}
{GRAPH ISOMORPHISM}
{Two graphs, $G$ and $G'$.}
{Are $G$ and $G'$ isomorphic?}
\end{problem}

GRAPH ISOMORPHISM is thought to be a good candidate
for a problem in \tup{NP} that is neither \npc\ nor
in \tup{P} (see~\cite{garey}).

A decision problem that is polynomially equivalent to
GRAPH ISOMORPHISM is \emph{isomorphism-complete}.
In this section we show that the analogous matroid problem is
in general isomorphism-complete.

\begin{problem}
{MATROID ISOMORPHISM}
{$(M,\, \tsc{Input})$ and $(M',\, \tsc{Input})$.}
{Are $M$ and $M'$ isomorphic?}
\end{problem}

\begin{lemma}
\label{lem13}
\tup{MATROID ISOMORPHISM} is polynomially reducible to
\tup{GRAPH ISOMORPHISM}\@.
\end{lemma}

\begin{proof}
A proof can be found in~\cite{mayhew2}, we give here an outline.
We must construct a polynomial-time computable transformation that
takes descriptions of matroids to graphs in such a way that
isomorphism is preserved.
There are many ways in which this can be accomplished.
The key idea is that a list of characteristic vectors, representing
subsets of the ground set, can be seen as the rows of the
vertex-adjacency matrix of a bipartite graph.

The rest of the demonstration involves refining the
transformation so that, given the unlabelled bipartite graph,
it is possible to reconstruct the matroid description,
up to relabelling.
This guarantees that the transformation preserves isomorphism.
Thus we must somehow distinguish the vertices that correspond to
subsets of the ground set from the vertices that correspond to
elements of the ground set.
In the case that \tsc{Input} relies upon ranks being
assigned to sets, as is the case when
$\tsc{Input} = \tsc{Cyclic Flats}$, we must find a method of
encoding binary representations of integers in the form of graphs.
Constructing a transformation that satisfies these criteria, and
confirming that it is polynomial-time computable, is an easy
exercise.
\end{proof}

Next we develop a polynomial transformation from graphs to
matroid descriptions.
Suppose that $G$ is a simple graph with $n$ vertices,
$\{v_{1},\ldots, v_{n}\}$, and $m$ edges,
$\{e_{1},\ldots, e_{m}\}$.
Assume that $n \geq 3$.
Let $X = \{x_{1},\ldots, x_{n}\}$,
$X' = \{x'_{1},\ldots, x'_{n}\}$, and
$Y = \{y_{1},\ldots, y_{m}\}$ be disjoint sets.
The matroid $\Phi(G)$ has rank $3$, and the ground set of $\Phi(G)$ is
$X \cup X' \cup Y$. The non-spanning circuits of $\Phi(G)$ are
sets of the form $\{x_{i},\, x'_{i}\}$, where
$i \in \{1,\ldots, n\}$, and the sets
\begin{displaymath}
\{\{z_{i},\, z_{j},\, y_{k}\} \mid z_{i} \in \{x_{i},\, x'_{i}\},\,
z_{j} \in \{x_{j},\, x'_{j}\},\ \mbox{and}\ e_{k}\ \mbox{joins}\ 
v_{i}\ \mbox{to}\ v_{j}\}.
\end{displaymath}
Thus $\Phi(G)$ can be formed by placing $n$ parallel pairs,
corresponding to the $n$ vertices of $G$, in the plane in
general position, and placing an element between parallel
pairs that correspond to adjacent vertices, in such a way
that no additional dependencies are formed.

Let $G$ be a graph.
The \emph{cyclomatic number of $G$} is $|V(G)| - |E(F)|$,
where $F$ is a spanning forest of $G$.
The \emph{bicircular matroid of a graph, $G$}, denoted by
$B(G)$, has the edge set of $G$ as its ground set.
The circuits of $B(G)$ are exactly the minimal connected
edge sets of $G$ with cyclomatic number two, known as
\emph{bicycles}.
Thus a set of edges is independent in $B(G)$ if and only
if the subgraph of $G$ it induces contains at
most one cycle in every connected component.

Suppose that $G$ is simple and has at least three vertices.
Let $G^{\circ\circ}$ be the graph which is obtained by adding two
loops at each vertex of $G$.
Then $\Phi(G)$ is isomorphic to the matroid obtained by
repeatedly truncating $B(G^{\circ\circ})$ so that its rank is
reduced to three.

Suppose that $G$ is a simple graph with $n \geq 3$ vertices and
$m$ edges.
Then $|\Phi(G)| = 2n + m$.
Since $m \leq n^{2}$ the number of independent sets in $\Phi(G)$
is $O(n^{6})$.
Thus $(\Phi(G),\, \tsc{Input})$ can be constructed in polynomial
time from a description of $G$ when
$\tsc{Independent Sets} \leq \tsc{Input}$.

It is easy to demonstrate that if $G$ and $G'$ are simple graphs
on at least $3$ vertices, then $G$ and $G'$ are isomorphic if and
only if $\Phi(G)$ and $\Phi(G')$ are isomorphic.
Since GRAPH ISOMORPHISM is generally defined in terms of simple
graphs we have proved the following result.

\begin{lemma}
\label{lem12}
If $\tsc{Independent Sets} \leq \tsc{Input}$, then
\tup{GRAPH ISOMORPHISM} is polynomially reducible to
\tup{MATROID ISOMORPHISM}\@.
\end{lemma}

Using duality we see that
$(\Phi(G)^{*},\, \tsc{Spanning Sets})$ can also
be constructed in polynomial time from a description
of $G$, where $G$ is any simple graph with
at least three vertices.
Since $G$ and $G'$ are isomorphic graphs if and only if
$(\Phi(G))^{*}$ and $(\Phi(G'))^{*}$ are isomorphic matroids
we have proved the following.

\begin{theorem}
\label{thm2}
If $\tsc{Spanning Sets} \leq \tsc{Input}$, or
$\tsc{Independent Sets} \leq \tsc{Input}$, then
\tup{MATROID ISOMORPHISM} is isomorphism-complete.
\end{theorem}

\section{Detecting minors}
\label{sct4}

It is a well known observation of Seymour's that
an oracle Turing machine cannot decide whether a matroid has a
$U_{2,4}$\dash minor in time that is polynomially-bounded
by the size of the ground set~\cite{seymour}.
We will show that, in general, deciding whether a matroid contains
a minor isomorphic to some fixed matroid can be done in polynomial
time, given a list of the independent sets or some similar input.
It is routine to verify the following result.

\begin{proposition}
\label{prop3}
Suppose that $M$ is a matroid and that $X$ and $Y$ are disjoint
subsets of $E(M)$.
It is possible to construct
$(M / X \del Y,\, \tsc{Input})$ in polynomial time given
$(M,\, \tsc{Input})$.
\end{proposition}

\begin{proposition}
\label{prop1}
Let $N$ be a matroid and let $t$ be the number of
circuits of $N$, where $t > 0$. If $M$ is a matroid on
the ground set $E$, and $M$ has a minor isomorphic to $N$
on the set $A \subseteq E$, then there exist
circuits, $C_{1},\ldots, C_{t}$, of $M$, such
that, if $X = (C_{1}\cup \cdots \cup C_{t}) - A$, and
$Y = E - (A \cup X)$, then $M/ X \del Y \iso N$.
\end{proposition}

\begin{proof}
Suppose that $M / X' \del Y'$ is isomorphic to $N$, where
$(X',\, Y')$ is a partition of $E - A$.
We may assume that $X'$ is independent.
There are exactly $t$ circuits, $C'_{1},\ldots, C'_{t}$,
in $M / X' \del Y'$.
For each circuit, $C'_{i}$, there exists a circuit $C_{i}$
of $M$ such that $C_{i} \subseteq C'_{i} \cup X'$ and
$C'_{i} = C_{i} - X'$.
Let $X$ be the set $(C_{1} \cup \cdots \cup C_{t}) - A$ and
let $Y$ be $E - (A \cup X)$.
It remains to show that $M / X \del Y \iso N$.
Note that $X \subseteq X'$.
If $X = X'$ then we are done so suppose that $x$ is an element
in $X' - X$.
Let $M' = M / (X' - x) \del (Y' \cup x)$.
To complete the proof it will suffice to show that
$M' = M / X' \del Y'$.
This is elementary.
\end{proof}

Let $N$ be some fixed matroid. It is a consequence of
Propositions~\ref{prop3} and~\ref{prop1}
that the following problem is, in general, in \tup{P}.

\begin{problem}
{DETECTING AN $N$\dash MINOR}
{$(M,\, \tsc{Input})$.}
{Does $M$ have a minor isomorphic to $N$?}
\end{problem}

\begin{theorem}
\label{thm3}
If $\tsc{Input} \leq \tsc{Circuits}$ or
$\tsc{Input} \leq \tsc{Hyperplanes}$, then
\tup{DETECTING AN $N$\dash MINOR} is in \tup{P} for any
fixed matroid $N$.
\end{theorem}

\begin{proof}
By duality it will suffice to prove the theorem only when
\tsc{Input} is equal to \tsc{Circuits}.
We may assume that we have the list of circuits of $N$,
for we can construct it in constant time.

Suppose that the ground set of $M$ is $E$, where $|E| = n$,
and that $M$ contains $c$ circuits, so that
$|(M,\, \tsc{Circuits})| = \Theta(nc)$. Suppose also that
$|E(N)| = s$, and that $N$ has exactly $t$ circuits. We may assume
that $t > 0$, for otherwise the problem is trivially in \tup{P}\@.

An algorithm to check whether $M$ has an $N$\dash minor
could simply work its way through every subset, $A$, of $E$,
such that $|A| = s$, and every set, $\{C_{1},\ldots, C_{t}\}$,
of $t$ circuits of $M$.
By Proposition~\ref{prop3} it is possible to construct
$(M / X \del Y,\, \tsc{Circuits})$, where
$X = (C_{1}\cup \cdots \cup C_{t}) - A$ and $Y = E - (A \cup X)$, and
then check whether $M / X \del Y \iso N$ in polynomial time.
Proposition~\ref{prop1} guarantees that if $M$ does have an minor
isomorphic to $N$, then this procedure will find such a minor.
Checking whether an isomorphism exists between $N$ and $M / X \del Y$
will take some constant time, so the running time of the algorithm is
determined by the number of subsets, $A \subseteq E$, and the number
of families of circuits we must check.
The first quantity is $\binom{n}{s}$, and the second is $\binom{c}{t}$,
so the running time of the algorithm is $O(n^{s}c^{t})$.
Since $s$ and $t$ are fixed constants the algorithm runs in polynomial time.
\end{proof}

In contrast to Theorem~\ref{thm3}, Hlin\v{e}n\'{y} shows that
problem of deciding whether $M$ has a minor isomorphic $N$ for a fixed
matroid $N$ is \npc\ when the input consists of a representation of
$M$ over the rational numbers~\cite{hlineny}.
On the other hand, Geelen, Gerards, and Whittle conjecture that
the problem is in \tup{P} when the input consists of a representation
of $M$ over a finite field~\cite{geelen3}.

The exponent in the running time of the algorithm described in
Theorem~\ref{thm3} depends upon $N$.
It is natural to ask whether there is a
\emph{fixed-parameter tractable} algorithm for
DETECTING AN $N$\dash MINOR, that is, an algorithm which runs
in time $f(N)|(M,\, \tsc{Input})|^{k}$, where $f$ is a function
which depends only on $N$ and $k$ is a fixed constant.
Certainly such an algorithm exists when $\tsc{Input} = \tsc{Rank}$, for
we can simply consider minors of the form $M / X \del Y$ for all
disjoint sets $X$ and $Y$, and check by brute force
whether any such minor is isomorphic to $N$.
This leads to an FPT algorithm.
The existence of such an algorithm when $M$ is described via an
input that lies above \tsc{Rank} is an open problem.

We have shown that detecting a fixed minor can be solved in polynomial
time, in general.
However, the problem of detecting a minor which forms part of the
input is in general \npc.

\begin{problem}
{MINOR ISOMORPHISM}
{$(M,\, \tsc{Input})$ and $(N,\, \tsc{Input})$.}
{Does $M$ have a minor isomorphic to $N$?}
\end{problem}

\begin{theorem}
\label{thm4}
If $\tsc{Independent Sets} \leq \tsc{Input}$
or $\tsc{Spanning Sets} \leq \tsc{Input}$, then
\tup{MINOR ISOMORPHISM} is \npc.
\end{theorem}

\begin{proof}
It is easy to see that MINOR ISOMORPHISM is in \tup{NP}.
The following problem is well known to be \npc.

\begin{problem}
{SUBGRAPH ISOMORPHISM}
{Two graphs, $G$ and $H$.}
{Does $G$ have a subgraph isomorphic to $H$?}
\end{problem}

Suppose that the graphs $G$ and $H$ correspond to an instance
of SUBGRAPH ISOMORPHISM\@.
We may assume that $G$ and $H$ are simple and that both have at
least three vertices.
Since we can construct either $(\Phi(G),\, \tsc{Input})$ or
$(\Phi(G)^{*},\, \tsc{Input})$ in polynomial time from the
description of $G$ (and the same statement applies for $H$),
the proof will be complete if we can demonstrate that $G$ contains
a subgraph isomorphic to $H$ if and only if $\Phi(G)$ contains a
minor isomorphic to $\Phi(H)$.
This is easily done.
\end{proof}

Suppose that \mcal{M} is a family of matroids.
A natural computational problem is to ask whether a matroid
$M$ has a minor isomorphic to a member of \mcal{M} with
specified size.
The proof of Theorem~\ref{thm4} can be used to show that
if $\mcal{M} = \{\Phi(K_{n}) \mid n \geq 1\}$ then this
problem is in general \npc.
We will conclude by showing that the problem is in general also \npc\
when $\mcal{M} = \{U_{r,n} \mid n \geq r\}$, where $r$ is a fixed
constant.

\begin{problem}
{$U_{r,n}$\dash MINOR}
{$(M,\, \tsc{Input})$ and an integer $n$.}
{Does $M$ have a minor isomorphic to $U_{r,n}$?}
\end{problem}

\begin{theorem}
\label{thm6}
If $r > 2$ is a fixed integer and
$\tsc{Independent Sets} \leq \tsc{Input}$ then
\tup{$U_{r,n}$\dash MINOR} is \npc.
\end{theorem}

\begin{proof}
It is easy to see that $U_{r,n}$\dash MINOR is in \tup{NP}\@.
We construct a reduction from the following \npc\ problem.

\begin{problem}
{INDEPENDENT SET}
{An integer $k$ and a graph $G$.}
{Does $G$ contain an independent set of $k$ vertices?}
\end{problem}

Let $r > 2$ be a fixed integer.
Suppose that the integer $k$ and the simple graph $G$ are an
instance of INDEPENDENT SET\@.
Let $m$ be the number of edges in $G$.
Let $t$ be $\lceil (r-1)/2 \rceil$, and let $tG^{\circ}$ be
the graph obtained by adding a loop to each vertex of $G$ and
then replacing each non-loop edge with a path of length $t$.
The matroid $\Phi_{r}(G)$ is the bicircular matroid of
$tG^{\circ}$, repeatedly truncated so that its rank is equal
to $r$.

Checking whether a set of edges of $tG^{\circ}$ is
independent in $B(tG^{\circ})$ can certainly be done in
polynomial time.
Since $r$ is a fixed integer, and no independent set of
$\Phi_{r}(G)$ exceeds $r$ in size, it follows that
$(\Phi_{r}(G),\, \tsc{Input})$ can be constructed in
polynomial time as long as $\tsc{Independent Sets} \leq \tsc{Input}$.

We complete the proof by showing that $G$ has an
independent set of $k$ vertices if and only
if $\Phi_{r}(G)$ has a rank\dash $r$ uniform minor
of size $k + mt$.

Every non-loop cycle of $tG^{\circ}$
contains at least $3t$ edges, so if a bicycle of
$tG^{\circ}$ contains at most one loop, then it has
size at least $3t + 1$.
This quantity is greater than $r$, as $r \geq 3$.
If a bicycle contains two loops, then either it contains
exactly $t + 2$ elements, or its size is at least $2t + 2$.
It is straightforward to confirm that $t + 2 \leq r$ and
that $2t + 2 > r$.
This shows that the non-spanning circuits of $\Phi_{r}(G)$
are exactly the sets containing the $t$ edges in a path joining
two vertices of $G$ and the two loops incident with those vertices.

Suppose that $G$ contains an independent set of
$k$ vertices.
These vertices correspond to $k$ loops of $tG^{\circ}$.
The set containing these loops and all the non-loop edges
of $tG^{\circ}$ cannot contain a non-spanning circuit of
$\Phi_{r}(G)$, by the above discussion.
Thus restricting $\Phi_{r}(G)$ to this set of $k + mt$ elements
gives a uniform minor.
We may assume that $m \geq 3$, so $k + mt$ is certainly no less
than $r$.
Therefore $\Phi_{r}(G)$ contains a rank\dash $r$ uniform minor
with $k + mt$ elements.

For the converse suppose that $\Phi_{r}(G)$ has a
rank\dash $r$ uniform minor on the set $A$, where
$|A| = k + mt$, and assume that $G$ has no
independent set of $k$ vertices.
The number of non-loop edges in $A$ is at most $mt$.
Suppose that $A$ has been chosen so that it contains as many
non-loop edges as possible.
Now $A$ contains at least $k$ loops, so there must be a pair of
loops, $l$ and $l'$, in $A$ that correspond to adjacent vertices in $G$.
Therefore one of the $t$ edges that join $l$ to $l'$ in $tG^{\circ}$
does not belong to $A$.
Let us call this edge $e$.
Then $(A - l) \cup e$ contains no non-spanning circuits of
$\Phi_{r}(G)$, and our assumption on $A$ is contradicted.
Therefore $G$ has an independent set of $k$ vertices. 
\end{proof}

If $r \leq 2$ then $U_{r,n}$\dash MINOR is trivially in \tup{P}\@.
Using duality, we can show that the problem of deciding whether
$M$ has a minor isomorphic to $U_{n-r,n}$ is \npc\ for fixed
values of $r > 2$ as long as $\tsc{Spanning Sets} \leq \tsc{Input}$.

\section{Open Problems}
\label{sct8}

In this summary section we collect some open problems.
$3$\dash MATROID INTERSECTION is known to be either in
\tup{P} or \npc\ for all but two of the types of input
mentioned in Section~\ref{sct2}.
The status of the problem is open when \tsc{Input}
is either \tsc{Cyclic Flats} or \tsc{Flats}.
MATROID ISOMORPHISM is \npc\ for all forms of input,
except possibly \tsc{Rank}.
Deciding if MATROID ISOMORPHISM can be solved in polynomial
time when the matroids are described via the rank of
each of their subsets is an open problem.
More generally, it would be interesting to know if there is
any `natural' computational problem which is
\npc\ for the \tsc{Rank} input.

DETECTING AN $N$\dash MINOR is known to be in \tup{P} for
all forms of input except \tsc{Non-Spanning Circuits},
\tsc{Cyclic Flats}, and \tsc{Dependent Hyperplanes}.
The status of the problem for these types of input
is unknown.
The existence or otherwise of an FPT algorithm for
MINOR ISOMORPHISM is known only when
$\tsc{Input} = \tsc{Rank}$.
Otherwise the problem is open.

\section{Acknowledgements}
\label{sct7}

I thank my supervisor, Professor Dominic Welsh, and
Professors Colin McDiarmid and James Oxley for very
useful advice and discussion.


\begin{thebibliography}{10}

\bibitem{edmonds3}
J.~Edmonds.
\newblock Submodular functions, matroids, and certain polyhedra.
\newblock In \emph{Combinatorial Structures and their Applications (Proc.
  Calgary Internat. Conf., Calgary, Alta., 1969)}, pp. 69--87. Gordon and
  Breach, New York (1970).

\bibitem{garey}
M.~R. Garey and D.~S. Johnson.
\newblock \emph{Computers and intractability : a guide to the theory of
  \textup{NP}-completeness}.
\newblock W. H. Freeman and Co., San Francisco, Calif. (1979).

\bibitem{geelen3}
J.~Geelen, B.~Gerards, and G.~Whittle.
\newblock Towards a structure theory for matrices and matroids.
\newblock In \emph{International Congress of Mathematicians. Vol. III}, pp.
  827--842. Eur. Math. Soc., Z\"urich (2006).

\bibitem{hausmann2}
D.~Hausmann and B.~Korte.
\newblock Oracle algorithms for fixed-point problems---an axiomatic approach.
\newblock In \emph{Optimization and operations research (Proc. Workshop, Univ.
  Bonn, Bonn, 1977)}, volume 157 of \emph{Lecture Notes in Econom. and Math.
  Systems}, pp. 137--156. Springer, Berlin (1978).

\bibitem{hausmann3}
D.~Hausmann and B.~Korte.
\newblock Worst-case behaviour of polynomial bounded algorithms for
  independence systems.
\newblock \emph{Z. Angew. Math. Mech.} \textbf{58} (1978), no.~7, T477--T479.

\bibitem{hausmann}
D.~Hausmann and B.~Korte.
\newblock Algorithmic versus axiomatic definitions of matroids.
\newblock \emph{Math. Programming Stud.}  (1981), no.~14, 98--111.

\bibitem{hlineny}
P.~Hlin{\v{e}}n{\'y}.
\newblock On matroid representability and minor problems.
\newblock In \emph{Mathematical foundations of computer science 2006}, volume
  4162 of \emph{Lecture Notes in Comput. Sci.}, pp. 505--516. Springer, Berlin
  (2006).

\bibitem{knuth}
D.~E. Knuth.
\newblock The asymptotic number of geometries.
\newblock \emph{J. Combin. Theory Ser. A} \textbf{16} (1974), 398--400.

\bibitem{lawler2}
E.~L. Lawler.
\newblock Polynomial-bounded and (apparently) non-polynomial-bounded matroid
  computations.
\newblock In \emph{Combinatorial algorithms (Courant Comput. Sci. Sympos. 9,
  New York Univ., New York, 1972)}, pp. 49--57. Algorithmics Press, New York
  (1973).

\bibitem{mayhew2}
D.~Mayhew.
\newblock \emph{Matroids and complexity}.
\newblock {DP}hil {T}hesis, {U}niversity of {O}xford (2005).

\bibitem{oxley}
J.~G. Oxley.
\newblock \emph{Matroid theory}.
\newblock Oxford University Press, New York (1992).

\bibitem{robinson}
G.~C. Robinson and D.~J.~A. Welsh.
\newblock The computational complexity of matroid properties.
\newblock \emph{Math. Proc. Cambridge Philos. Soc.} \textbf{87} (1980), no.~1,
  29--45.

\bibitem{seymour}
P.~D. Seymour.
\newblock Recognizing graphic matroids.
\newblock \emph{Combinatorica} \textbf{1} (1981), no.~1, 75--78.

\end{thebibliography}
\end{document}